\documentclass[11pt]{article}
\usepackage[utf8]{inputenc}
\usepackage{amsmath}
\usepackage{geometry}
\usepackage{graphicx}
\usepackage{anyfontsize}
\usepackage{algorithm}
\usepackage{algorithmic}
\usepackage{tikz}
\usetikzlibrary{calc}
\usepackage{pgfplots}
\usepackage{graphicx}
\usepackage{ amssymb }
\usepackage{subcaption}
\selectfont

\title{Exploring the Ant Mill: Numerical and Analytical Investigations of Mixed Memory-Reinforcement Systems}
\author{Ria Das \\ Phillips Exeter Academy \\ \\ Mentored by Andrew Rzeznik, MIT \\ \\ MIT PRIMES}	

\date{}

\begin{document}
\maketitle
\begin{abstract}
Under certain circumstances, a swarm of a species of trail-laying ants known as army ants can become caught in a doomed revolving motion known as the death spiral, in which each ant follows the one in front of it in a never-ending loop until they all drop dead from exhaustion. This phenomenon, as well as the ordinary motions of many ant species and certain slime molds, can be modeled using reinforced random walks and random walks with memory. In a reinforced random walk, the path taken by a moving particle is influenced by the previous paths taken by other particles. In a random walk with memory, a particle is more likely to continue along its line of motion than change its direction. Both memory and reinforcement have been studied independently in random walks with interesting results. However, real biological motion is a result of a combination of both memory and reinforcement. In this paper, we construct a continuous random walk model based on diffusion-advection partial differential equations that combine memory and reinforcement. We find an axi-symmetric, time-independent solution to the equations that resembles the death spiral. Finally, we prove numerically that the obtained steady-state solution is stable.
\end{abstract}

\section{Introduction}

The process by which living organisms respond to chemical stimuli in their environments by moving towards or away from them is called \textit{chemotaxis}. Movement towards a higher concentration of chemical is known as positive taxis (\textit{taxis} in Greek means \textit{to arrange}), while movement towards a lower concentration of chemical is known as negative taxis (Friedman and Tello, 2002). Chemotaxis has been observed in several species of bacteria, including the well-studied soil species known as the myxobacteria; slime molds; and ant species such as army ants (Othmer and Stevens, 1997). A chemotaxis process also occurs in the growth of new blood vessels, known as \textit{angiogenesis}, that occurs in the expansion of a tumor. In tumor-induced angiogenesis, chemical agents secreted by a tumor attract neighboring endothelial cells, which form the surface of blood vessels that then sprout towards the tumor to provide nourishment (Friedman and Tello, 2002). \\

The mathematics that underlies each of these chemotaxis scenarios is based on an object known as the random walk. A random walk is any path that consists of a sequence of random steps. Examples of random walks include the path traveled by a molecule in a liquid or gas and the path of a foraging animal. The class of random walks that have been used in the past to model chemotaxis processes are called reinforced random walks (RRWs). In this type of random walk, the walker can modify (reinforce) the favorability of traveling along a particular path, for example by dropping an attracting chemical substance. A variety of different reinforced random walk systems based on differing reinforcement mechanisms have modeled many observed chemotaxis situations with success (Pemantle, 2006). \\

Although reinforced random walk models have been able to replicate several different chemotaxis phenomena, there remain limitations to these systems (Codling, Plank, and Benhamou, 2008). In particular, reinforced random walks do not take into account another important property of animal and cell motion: the tendency of these bodies to continue moving in the same direction. This property of inertia is explicitly considered in a class of random walks called correlated random walks (CRWs). In these walks, there is a correlation between successive step orientations, which produces a local directional bias that causes each step to tend to point in the same direction as the previous one. The influence of the initial direction of motion gradually diminishes over time. \\

\begin{figure}
\begin{center}
\includegraphics[scale=0.15]{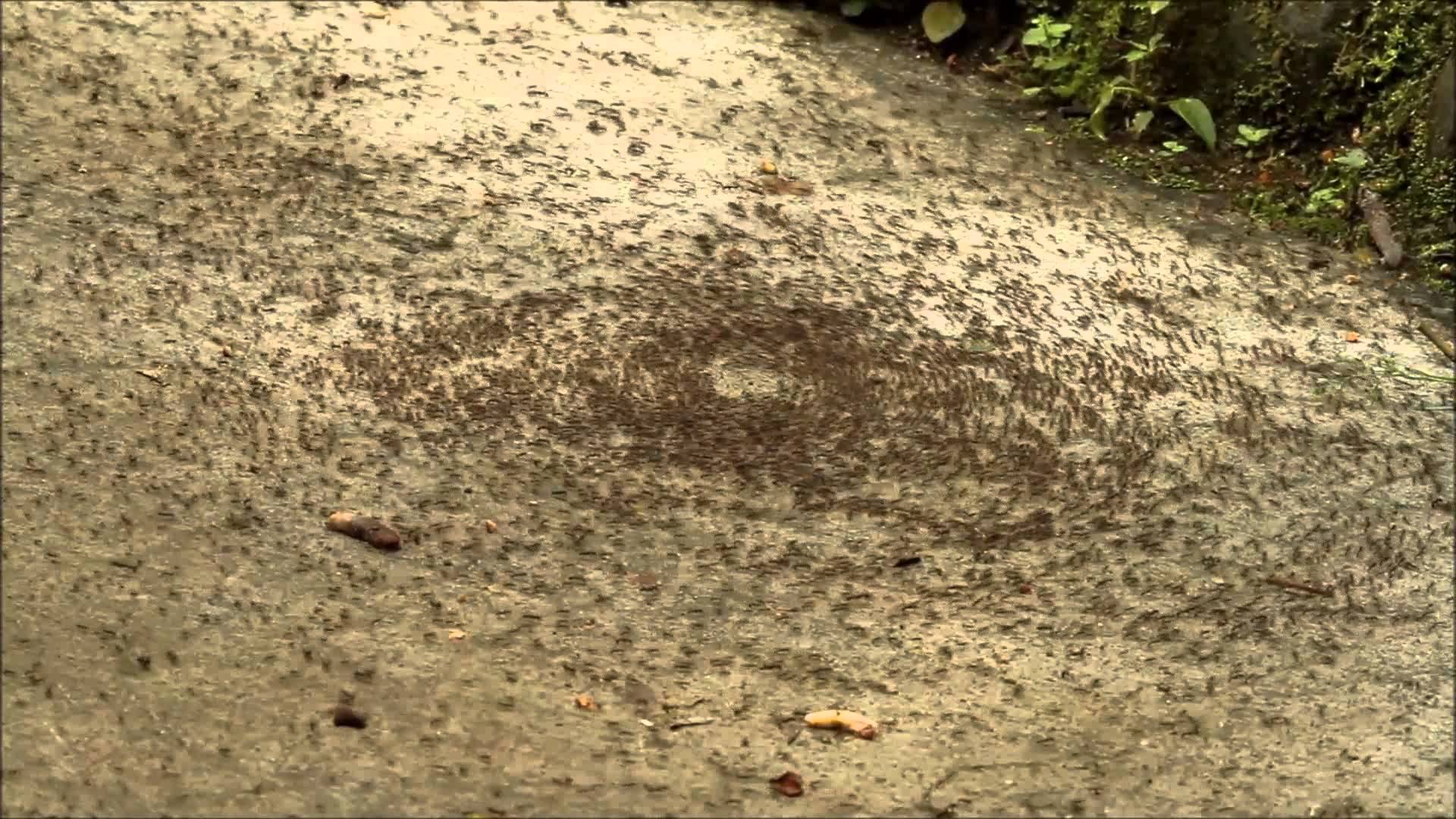}
\end{center}
\caption{An ant death spiral. [1]}
\end{figure}

\begin{figure}[h!]
\begin{center}
\begin{tikzpicture}[fill opacity=.6,draw opacity=1, fillednode/.style={circle, minimum       size=1pt, inner sep=1pt, fill},blanknode/.style = {circle, minimum size=0pt,inner sep=0pt, fill}]
	\node[fillednode] (A) at (0,0) [label=north east:] {};	
	\draw[blue,thick,dashed] (0,0) circle (2cm);
	\shade[inner color=black,outer color=white] (0,0) circle (3cm);
	\draw[red, very thick] (0,-2) circle (0.03cm);
	\draw[white, very thick] (0,0) circle (0.1cm);
	\draw[white, very thick] (0,0) circle (0.13cm);
	\draw[white, very thick] (0,0) circle (0.16cm);
	\draw[white, very thick] (0,0) circle (0.19cm);
	\draw[white, very thick] (0,0) circle (0.22cm);
	\draw[white, very thick] (0,0) circle (0.25cm);
	\draw[white, very thick] (0,0) circle (0.28cm);
	\draw[white, very thick] (0,0) circle (0.31cm);	
	\draw[white, very thick] (0,0) circle (0.34cm);		
	\draw[red,thick,->] (0,-2) -- (0,-0.8);
	\draw[red,thick,->] (0,-2.01) -- (1.2,-2.01);
\end{tikzpicture}
\caption{Schematic of ant death spiral. Gray background shading indicates concentration of chemical---darker shading indicates higher concentration and lighter shading indicates lower concentration. Centripetal and tangential influences on a single ant's motion are shown by red arrows.}
      \label{fig:memIntersection4}
\end{center}
\end{figure}
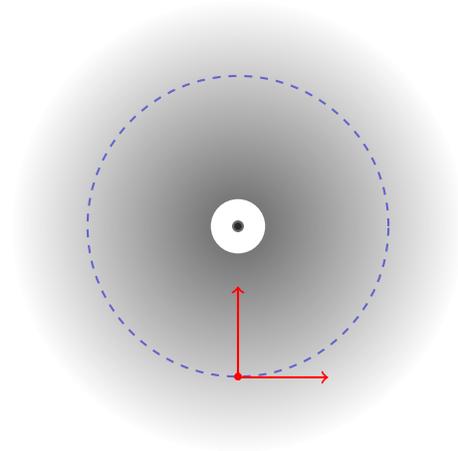

The motion of real biological systems is complex and is based on both correlation, which is also termed memory, and reinforcement. The fact that these two factors influencing motion have only been studied independently may explain why there continue to be observed phenomena that have not been replicated by past models. One such phenomenon is the army ant death spiral. In this example, we suspect that there are two basic directional forces influencing ant motion: $(1)$ a chemical-gradient-induced centripetal force and $(2)$ a tangential force related to inertia. These two influences are shown for a single ant caught in the spiral in Figure \ref{fig:memIntersection4}. The concentration of chemical is highest in a ring around the center of the spiral and fades out radially, creating a gradient that points toward the center for every point that the ants are located. The ``eye" of the spiral corresponds to the region where no ants are located in Figure $1$. The combined effect of this centripetal force and the tangential pull of inertia produces the observed circular motion of every ant in the configuration. \\

If reinforcement based on the trails laid by the ants was the only factor influencing their motion, it is unlikely that they would maintain circular motion. Instead, they would all move closer to the center of the spiral, where the chemical concentration was higher. Similarly, if the chemical gradient had a negligible effect, and inertia was the only bias propelling the ants, they would also not have circular motion. We seek a mathematical description that finds a balance between these two forces so as to replicate the army ant spiral that we observe in real life. \\

In this paper, we develop mixed memory-reinforcement models with the aim of replicating the army ant death spiral. Moreover, we aim to prove the stability of the ant death spiral mathematically based on the developed models. We begin our exploration of mixed memory and reinforcement in random walks with intuitive models based on integral equations (Section $2$). Due to the determined limitations of these systems of equations, we then consider random walk models based on diffusion-advection equations. We construct a system of diffusion-advection partial differential equations that has a time-independent solution resembling the death spiral, and then numerically prove the solution is stable in time (Section $3$).

\section{Integral equation model}
We begin our quest for a mathematical model of the ant mill by developing an intuitive system that describes all properties of ant motion as simply as possible. In this model, we let the position $\vec{x}$ of an ant in a $2$-dimensional plane be independent of the ant's velocity, $\vec{v}$. For every $\vec{x}$ in position space, every $\vec{v}$ in velocity space, and every time $t>0$, the density $\rho(\vec{x},\vec{v},t)$ is the density of ants at that point in the position-velocity phase space at time $t$. The system thus keeps track of position and velocity information for every ant in its domain over time. We can simplify this model further by placing the additional restriction that all ants travel at a fixed speed, call it $v$. This assumption is reasonable because there is little variance in the cruising speeds of ants. Adopting it also allows us to define $\rho$ as a function simply of the angle $\theta$ of its velocity rather than of the vector $\vec{v}$, so that $\rho=\rho(\vec{x},\theta, t)$. \\ 

To define an equation for the evolution of $\rho(\vec{x},\theta,t)$, we must consider the factors that might change the density of particles at a particular location $\vec{x}$ and moving in the specific direction $\theta$ after an infinitesimal step in time. We draw this from the taxis model described by Codling, et. al. in their 2008 review paper on random walk models in biology. In their model, there exists some probability per time $\alpha$ that the velocity $\vec{v}$ of a particle will change at any time. Further, a \textit{reorientation kernel} $T(\vec{v},\vec{v}')$ is defined to be the probability that, given a reorientation occurs, the velocity of the particle will change from $\vec{v}'$ to $\vec{v}$. Since $T(\vec{v},\vec{v}')$ gives a probability, it satisfies the normalization condition that $\int T(\vec{v},\vec{v}')=1$. The evolution equation arising from these definitions in the Codling paper is
\begin{equation}
\frac{d\rho}{dt}=\frac{\partial \rho}{\partial t}+\vec{v}\cdot\nabla_{x}\rho = -\alpha\rho+\alpha\int T(\vec{v},\vec{v}')\rho(\vec{x},\vec{v}',t) dv',
\end{equation}

\noindent where $\nabla_x$ denotes the spatial gradient operator. Converting to our form of the system with fixed speed $v$ and $\rho(\vec{x},\theta,t)$ instead of $\rho(\vec{x},\vec{v},t)$, we have
\begin{equation}
\frac{d\rho}{dt}=\frac{\partial \rho}{\partial t}+\vec{v}\cdot\nabla_{x}\rho = -\alpha\rho+\alpha\int_{-\pi}^{\pi} T(\theta,\theta')\rho(\vec{x},\theta',t)d\theta'.
\end{equation}

In this equation, the left-hand side is the total derivative of 
$\rho$ with respect to time, equal to the partial derivative of $\rho$ in time plus the convective derivative ($\vec{v}\cdot\nabla_{x}\rho$) in the center expression. The right-hand side of the equation defines the mechanism of change in the particle density $\rho$: The first term on the right side, $-\alpha\rho$, signifies that $\alpha\rho$ is the rate at which particles at position $\vec{x}$ and traveling in direction $\theta$ randomly change direction, and hence are no longer part of the $\rho(\vec{x},\theta,t)$ phase density. The second integral term can be understood as the rate at which particles traveling in directions $\theta' \neq \theta$ randomly change their direction to $\theta$. Thus, we have a simple form for an evolution equation for the particle density $\rho$ in this definition. \\

The integration of both memory and reinforcement into our system, which we believe will produce the spiraling solution of the ant mill, is made in the specific form of the reorientation kernel $T$. In the Codling paper, the kernel $T$ is defined only in terms of the original direction $\theta'$ and the new direction $\theta$. Clearly, this setup serves as a strong model for a correlated random walk, as the probability of a particle engaging in that direction change could be made higher for small values of $|\theta-\theta'|$ and lower for large values of $|\theta-\theta'|$. We seek, however, a kernel $T$ that combines correlation with reinforcement based on the surrounding concentrations of chemical. In order to construct such a kernel, we must first define an evolution equation for the chemical concentration. We give a simple chemical evolution equation below, where $g(\vec{x},t)$ denotes the concentration:
\begin{equation}
\frac{\partial g}{\partial t}=\gamma\int_{-\pi}^{\pi}\rho(x,\theta',t)d\theta'-\beta g.  
\end{equation}
In this equation, in which $\gamma$ and $\beta$ are constants, the first integral term indicates the rate at which the particles located at $\vec{x}$ deposit chemical, and the second term signifies the rate at which the chemical evaporates. In essence, the equation describes that the rate at which the chemical concentration is changing at a location $\vec{x}$ at a time $t$ is equal to rate at which ants are adding chemical minus that at which the chemical is evaporating (an exponential rate). Given this definition, we introduce the following form for the reorientation kernel $T$:

\begin{equation}
T(\theta,\theta_{g})=\frac{1}{4\pi^2}(Jcos(\theta-\theta_{g})+1),
\end{equation}

\noindent where $\theta_{g}$ is the angle of the $\nabla g$ vector and $J$ is a constant. The reasoning behind this choice of kernel is as follows: The closer to $\theta-\theta_g$ is to $0$, the greater the value of $T(\theta,\theta_{g})$ due to the nature of the cosine function. This means that the closer the angle $\theta$, the intended level of direction, is to the gradient, the more likely it is for a particle initially in the $\theta'$ direction to change direction to $\theta$. Although we could also used $\theta'$ in our definition of $T$, using only $\theta$ and $\theta_{g}$ simplifies the kernel while still making sense physically, so we proceed with the current form. The constant $\frac{1}{4\pi^2}$ serves as a normalization factor, and $|J|<1$ so that $T>0$. \\

We hence have the following integral equation for $\rho$:
\begin{equation}
\frac{\partial \rho}{\partial t}+v(\cos\theta,\sin\theta)\cdot\nabla_{x}\rho=-\alpha\rho+\frac{\alpha}{4\pi^2}(Jcos(\theta-\theta_{g})+1)\int_{-\pi}^{\pi}\rho(x,\theta',t)d\theta'. 
\end{equation}
Our next step is to attempt to solve this equation for a steady-state, radially symmetric solution that resembles the ant mill. The symmetry of this solution form makes it valuable to convert it from rectangular to polar coordinates $(r,\phi)$. We notice that in the steady-state solution, the direction of the chemical gradient $\theta_{g}$ is always equal to the polar angle $\phi$, because the higher concentration of chemical in the center of the spiral pulls the ants around in a circle. Using these facts, letting $\phi=\theta_{g}=0$ without loss of generality, and simplifying, we arrive at the following equation from our original model,
\begin{equation}
(v\cdot2\pi i k \cos\theta +\alpha) P(ik,\theta)-\frac{\alpha}{4\pi^2}(Jcos\theta+1)\int_{-\pi}^{\pi} P(ik,\theta')d\theta'=0.
\end{equation}
\noindent where $P(ik,\theta)$ is the Fourier transform of $\rho$ (i.e. $P(ik,\theta)=\mathcal{F}(\rho)$). \\ %The details of this simplification process are given in Appendix A; we only give the highlights in this section for the sake of concision. 

The above integral equation is a homogeneous Fredholm equation of the second kind. According to Arkfen in Section $16.3$ of his book \textit{Mathematical Methods for Physicists}, however, this equation has no nontrivial solution: For all space, $0$ is a solution to the equation, so we cannot find a nonzero solution. This result means that the current formulation of our system of equations modeling ant motion has no steady-state. In hindsight, this was to have been expected---it was unlikely that a spiraling time-independent solution could have been obtained without a nonlinearity in the equations causing the ants to remain in orbit. \\

In addition, the inability of the integral equation model discussed in this section to produce a nontrivial steady-state solution suggests that the fundamental assumption underpinning the model---that the density $\rho$ is a function of both position $\vec{x}$ and velocity $\vec{v}$---is ill-fated. Indeed, storing a density value for every position and velocity in space is a lot of information, as we are basically keeping track of all of the individual particles that are moving around together in the domain. In the next section, we remedy this potential problem of too much information by developing a slightly simpler model. In our new model, we consider the average velocity of particles at every point in space, rather than all the velocities of the different particles at that point. Furthermore, the inability of our linear integral equation to produce a nontrivial solution tells us that nonlinearity plays a critical role in describing circular motion. Thus, this aspect of our model must also be fixed, by the addition of nonlinearity. Since nonlinear integral equations are difficult to work with, we turn to non-linear diffusion-advection equations to model our system instead, because they preserve much of the core physics while also being more pliable. In combination, these changes ultimately end up producing the steady-state result that we seek. \\

Based on what we have learned from the integral equation model, we now develop a new, nonlinear model based on diffusion-advection equations.

\section{Diffusion-advection model} \label{continuum}
\subsection{Construction}
\noindent To develop our new model of ant motion, we draw from past models in which the density $\rho$ is a function of only position and time, not velocity. The basis of our new system is the following pair of pure reinforcement equations for the particle density $\rho$ and chemical concentration $g$ in a one-dimensional environment from Othmer and Stevens (1997): 
\begin{equation}
\frac{\partial \rho}{\partial t} = D\frac{\partial}{\partial x}\left(\frac{\partial \rho}{\partial x}-\rho \frac{\beta}{\alpha+\beta g}\frac{\partial g}{\partial x}\right)
\end{equation}
\begin{equation}
\frac{\partial g}{\partial t} = \lambda\rho-g
\end{equation} \\
The first equation describes the evolution of the particle density and the second describes the evolution of the chemical concentration. The constant $\lambda$ indicates the magnitude of the rate of chemical deposition by the ants, and $D$ is the diffusion constant. We chose this particular system to serve as the basis for our new model because Othmer and Stevens have shown analytically and numerically that simple solutions of these equations are asymptotically stable. Further information on these simple solutions and their stability results can be found in Other and Stevens' 1997 paper. For now, we press onward to transform this pure-reinforcement system into a mixed memory-reinforcement one by adding a formulation of velocity. \\

In Othmer and Stevens' model, there is no concept of velocity, as indicated by the absence of a velocity variable $\vec{v}$ in the equations. Essentially, particles immediately forget their velocity at every instant, so there is no memory in the system. We change this by first adding a convective derivative to the left-hand side of the $\rho$ equation, so the left-hand side becomes the total derivative rather than the partial derivative. We also, of course, convert the one-dimensional equation into a two-dimensional one by changing the $\frac{\partial}{\partial x}$ operators to $\nabla$ operators:
\begin{equation}
\frac{d\rho}{dt} = \frac{\partial \rho}{\partial t} + \vec{v} \cdot \nabla\rho = \nabla\cdot(\nabla\rho -\rho \frac{\beta}{\alpha+\beta g} \nabla g). 
\end{equation}
In addition, we must define a velocity evolution equation to complete the memory-reinforcement system. It seems reasonable that the rate at which the velocity at a location $\vec{x}$ changes should be related to the concentration of chemical at that location. In particular, it would make sense that the derivative of the velocity $\vec{v}$ is in the same direction as the chemical gradient $\vec{g}$, since particles are attracted to higher concentrations of chemical, and that the steeper the gradient, the higher the rate of change. Taking these assumptions into account, we develop the simple velocity equation \\
\begin{equation}
\frac{d\vec{v}}{dt}=\frac{\partial \vec{v}}{\partial t}+\vec{v}\cdot\nabla \vec{v} = b\nabla g
\end{equation}

\noindent where $b$ is a constant. \\

Hence, we have constructed the following system of three partial differential equations that models memory-reinforcement ant motion from the original system by Othmer and Stevens:
\begin{equation} \label{eq:11}
\frac{d\rho}{dt} = \frac{\partial \rho}{\partial t} + \vec{v} \cdot \nabla\rho = \nabla\cdot(\nabla\rho -\rho \frac{\beta}{\alpha+\beta g} \nabla g). 
\end{equation} 
\begin{equation} \label{eq:12}
\frac{\partial g}{\partial t} = \lambda\rho-g
\end{equation}
\begin{equation} \label{eq:13}
\frac{d\vec{v}}{dt}=\frac{\partial \vec{v}}{\partial t}+\vec{v}\cdot\nabla \vec{v} = b\nabla g.
\end{equation}
Our next step is to solve our new system for a steady-state solution that resembles the ant mill.

\subsection{Steady-state solution}
\noindent We seek a solution to our system that resembles the ant death spiral. Such a solution must be radially symmetric; in other words, the density $\rho$, the chemical concentration $g$, and the velocity $\vec{v}$ must have no $\theta$ dependence. Since each of $\rho$, $g$, and $\vec{v}$ must therefore only depend on the radial distance $r$, we can rewrite the system of equations \ref{eq:11}, \ref{eq:12}, and \ref{eq:13} as 
\begin{equation}
\frac{d\rho}{dt}= \frac{\partial}{\partial r}\left(r\frac{\partial\rho}{\partial r}-r\rho\frac{\beta}{\alpha+\beta g}\frac{\partial g}{\partial r}\right)
\end{equation}
\begin{equation}
\frac{\partial g}{\partial t} = \lambda\rho-g
\end{equation}
\begin{equation}
\frac{d\vec{v}}{dt}=\frac{\partial \vec{v}}{\partial t}-\frac{v_{\theta}^2}{r}=b \frac{\partial g}{\partial r} 
\end{equation}
\noindent using the polar forms of the gradient, divergence, and material derivative operators, and where $v_{\theta}$ denotes the $\theta$ component of $\vec{v}$. Further noting that $\frac{d\rho}{dt}=0$, $\frac{\partial g}{\partial t}=0$, and $\frac{\partial \vec{v}}{\partial t}=0$ in the time-independent case, we rewrite the equations as \\
\begin{equation} \label{eq:17}
0=\frac{\partial}{\partial r}\left(r\frac{\partial\rho}{\partial r}-r\rho\frac{\beta}{\alpha+\beta g}\frac{\partial g}{\partial r}\right)
\end{equation}
\begin{equation} \label{eq:18}
0=\lambda\rho-g
\end{equation}
\begin{equation} \label{eq:19}
-\frac{v_{\theta}^2}{r}=b \frac{\partial g}{\partial r}.
\end{equation}
\noindent We can rearrange equation \ref{eq:18} to get 
\begin{equation}
g=\lambda\rho
\end{equation}
and we substitute this expression for $g$ into equation \ref{eq:17} to obtain the single-unknown equation \\
\begin{equation} \label{eq:21}
0= \frac{\partial}{\partial r}\left(r\frac{\partial\rho}{\partial r}-r\lambda\rho\frac{\beta}{\alpha+\beta\lambda \rho}\frac{\partial \rho}{\partial r}\right).
\end{equation} \\
\noindent To solve equation \ref{eq:21} for $\rho$, we first integrate both sides with respect to $r$, so we obtain \\
\begin{equation*}
\int 0 dr = \int\frac{\partial}{\partial r}\left(r\frac{\partial\rho}{\partial r}-r\lambda\rho\frac{\beta}{\alpha+\beta\lambda \rho}\frac{\partial \rho}{\partial r}\right)dr
\end{equation*}
\begin{equation}
C_{1} = r\frac{\partial\rho}{\partial r}-r\lambda\rho\frac{\beta}{\alpha+\beta\lambda \rho}\frac{\partial \rho}{\partial r}
\end{equation}
\noindent where $C_{1}$ is a constant. We multiply both sides of the equation by $(\alpha+\beta\lambda\rho)$ and simplify as follows:
\begin{equation*}
C_{1}(\alpha+\beta\lambda \rho) = r\frac{\partial\rho}{\partial r}(\alpha+\beta\lambda \rho)-r\lambda\rho\beta\frac{\partial \rho}{\partial r}
\end{equation*}
\begin{equation*}
C_{1}(\alpha+\beta\lambda \rho) = r\alpha\frac{\partial\rho}{\partial r}
\end{equation*}
\begin{equation}
r\frac{\partial \rho}{\partial r}+C_{1}\frac{\beta}{\alpha}\lambda\rho = -C_{1} 
\end{equation}
\noindent where the last equation follows by changing the sign of the arbitrary constant $C_{1}$. This final form of the equation is straightforward to solve for $\rho$ using separation of variables. We subtract $C_{1}\frac{\beta}{\alpha}\lambda\rho$ from both sides of the equation, and then divide the right-hand side by $(-C_{1}-C_{1}\frac{\beta}{\alpha}\lambda\rho)$ and $\frac{r}{dr}$ so that each side of the equation is in terms of one variable:
\begin{equation*}
r\frac{\partial \rho}{\partial r}= -C_{1}-C_{1}\frac{\beta}{\alpha}\lambda\rho 
\end{equation*}
\begin{equation}
\frac{1}{-C_{1}-C_{1}\frac{\beta}{\alpha}\lambda\rho}\partial \rho= \frac{1}{r}\partial r
\end{equation}
\noindent We then integrate both sides and simplify to obtain
\begin{equation*}
\int \frac{1}{-C_{1}-C_{1}\frac{\beta}{\alpha}\lambda\rho}\partial \rho= \int \frac{1}{r}\partial r
\end{equation*}
\begin{equation*}
-\frac{1}{C_{1}}\displaystyle {\int} \frac{1}{1+\frac{\beta}{\alpha}\lambda\rho}\partial \rho= \int \frac{1}{r}\partial r
\end{equation*}
\begin{equation} \label{eq:25}
-\frac{1}{C_{1}}\frac{\alpha}{\beta\lambda}\ln(1+\frac{\beta}{\alpha}\lambda\rho)= \ln r+C_{2}
\end{equation}
\noindent where $C_{2}$ is a constant. Equation \ref{eq:25} can be rearranged to obtain the final form
\begin{equation}
\rho = \frac{\alpha}{\beta\lambda}\left(C_{2}r^{-C_{1}-\frac{\alpha}{\beta\lambda}}-1\right)
\end{equation}
\noindent where the constant $C_{2}$ is redefined as $e^{C_{1}C_{2}}$, where the $C_{2}$ in the exponent is the original $C_{2}$. \\

\noindent This solution for $\rho$ is reasonable because, although $\rho$ approaches infinity at $r=0$, the integral of $\rho$ converges for $r=1$ to $r=\infty$ for values of the exponent. This models the death spiral because it would not make sense for moving particles to be circling extremely fast in the center of a physical spiral. This is observed in real ant death spirals, in which there is a certain radius in the center of the spiral within which no ants appear to be circling (see the ant death spiral with the empty circular center in Figure $1$ and the schematic in Figure $2$). \\

\noindent From our solution for $\rho$, since $g=\lambda\rho$, we also have
\begin{equation}
g = \frac{\alpha}{\beta}\left(C_{2}r^{-C_{1}-\frac{\alpha}{\beta\lambda}}-1\right)
\end{equation}
\noindent as the solution for $g$. We then use this $g$ solution to obtain the solution for $\vec{v}$ by manipulating equation \ref{eq:19} with a little algebra, as below: \\
\begin{equation*}
-\frac{v_{\theta}^2}{r}=b \frac{\partial g}{\partial r} 
\end{equation*}
\begin{equation*}
v_{\theta}=\sqrt{-br\frac{\partial g}{\partial r}}
\end{equation*}
\begin{equation}
v_{\theta}=\left(\sqrt{bC_{2}\left(C_{1}+\frac{\alpha}{\beta\lambda}\right)}\right)r^{-\frac{1}{2}\left(C_{1}+\frac{\alpha}{\beta\lambda}\right)}.
\end{equation}
\noindent Hence, the final steady-state, radially symmetric solution to the entire system is 
\begin{equation}
\rho = \frac{\alpha}{\beta\lambda}\left(C_{2}r^{-C_{1}-\frac{\alpha}{\beta\lambda}}-1\right)
\end{equation}
\begin{equation}
g = \frac{\alpha}{\beta}\left(C_{2}r^{-C_{1}-\frac{\alpha}{\beta\lambda}}-1\right)
\end{equation}
\begin{equation}
v_{\theta}=\left(\sqrt{bC_{2}\left(C_{1}+\frac{\alpha}{\beta\lambda}\right)}\right)r^{-\frac{1}{2}\left(C_{1}+\frac{\alpha}{\beta\lambda}\right)}.
\end{equation} \\
\noindent Since the solution approaches infinity as $r$ approaches $0$, we note that the radial domain must be $r>r_{a}$ for a positive $r_{a}$. For simplicity, when performing stability analysis, we limit the domain further by specifying an outer boundary $r_{b}$, so we have a finite domain of $[r_{a},r_{b}]$. \\

By considering the signs of the constants in this steady-state solution, we know that the solution must decay to $0$ as $r\rightarrow \infty$. We arrive at this conclusion by first noting that each of the constants $\alpha$, $\beta$, and $\lambda$ are positive, based on the nature of the phenomena we are modeling and, specifically, the fact that the ants are attracted to (not repelled by) the chemical. This means that the coefficients $\frac{\alpha}{\beta\lambda}$ and $\frac{\alpha}{\beta}$ in the $\rho$ and $g$ equations are positive, so the coefficient $C_{2}$ of $r$ must also be positive for the density and chemical concentration to be positive. Since $C_{2}$ is positive, $C_{1}+\frac{\alpha}{\beta\lambda}$ must also be positive, because otherwise the velocities $\vec{v}$ will be imaginary due to the square root in the coefficient $\sqrt{bC_{2}(C_{1}+\frac{\alpha}{\beta\lambda})}$ in the $\vec{v}$ equation. Finally, since the $C_{1}+\frac{\alpha}{\beta\lambda}$ is positive, the powers of $r$ in all three equations must be negative, indicating that the solutions all decay as $r\rightarrow \infty$. The fact that these coefficients confirm our expectation of what the solution should be is a promising sign of the robustness of the model. \\

Having arrived at the steady state, we can now perturb these time-independent solutions slightly and observe how the system behaves over time in response to determine its stability.

\subsection{Stability Analysis}

Before beginning stability analysis, we shall briefly review what stability analysis is and what it can show about a system. Essentially, in a time-dependent system of differential equations, a solution to the system can either blow up, collapse, or oscillate as time approaches infinity. Blowup refers to the scenario when the solution approaches positive or negative infinity as time increases. Collapse indicates when a solution approaches a finite value as time increases. In a system that oscillates as time approaches infinity, the limit of the solution does not exist (for $t\rightarrow \infty$). \\

\noindent To analyze the stability of our system of partial differential equations around the obtained steady-state solutions, we first let the steady-state solutions for $\rho$, $g$, and $\vec{v}$ be represented by $\rho_{0}(r)$, $g_{0}(r)$, and $\vec{v}_{0}(r)$. We define $\tilde{p}(r,\theta,t)$, $\tilde{g}(r,\theta,t)$, and $\tilde{v}(r,\theta,t)$ to be perturbations to the steady state equations, each of which is a function of $r$, $\theta$, and $t$. Each perturbation is of order $\epsilon_{1}>0$ where $\epsilon_{1}$ is a small real number. Thus we can write

\[\rho = \rho_{0}+\epsilon_{1}\tilde{\rho},\]
\begin{equation}
g = g_{0}+\epsilon_{1}\tilde{g},
\end{equation}
\[\vec{v} = \vec{v}_{0}+\epsilon_{1}\tilde{v}.\]

\noindent We seek expressions for the unknown perturbations $\tilde{p}$, $\tilde{g}$, and $\tilde{v}$ that solve the general (not time-independent and axi-symmetric) system. Analyzing these expressions will allow us to assess whether small perturbations to the steady-state solutions cause the system to blow up or stabilize in infinite time. To solve for $\tilde{p}$, $\tilde{g}$, and $\tilde{v}$, we will first linearize the time-dependent equations for $\rho$ and $\vec{v}$, both of which are nonlinear. We do this because it is far more straightforward to solve linear differential equations than non-linear differential equations, both analytically and numerically. The full details of the linearization process for the $\rho$ and $\vec{v}$ equations are given in Appendix A; in this section, we omit some of the algebra for the sake of concision. \\

First, we convert the time-dependent equation for $\rho$ to polar coordinates, a more convenient coordinate system to use for modeling the spiraling motion of ants:
\begin{equation}
\frac{\partial \rho}{\partial t} + \vec{v} \cdot \nabla\rho = \nabla\cdot(\nabla\rho -\rho \frac{\beta}{\alpha+\beta g} \nabla g, 
\end{equation}
\vspace*{-1.32cm}
\begin{center}
\begin{equation} \label{eq:42}
\begin{split}
\frac{\partial \rho}{\partial t} & + v_{r}\frac{\partial \rho}{\partial r}+\frac{v_{\theta}}{r}\frac{\partial \rho}{\partial \theta}=\frac{1}{r}\frac{\partial \rho}{\partial r}+\frac{\partial^2 \rho}{\partial r^2}+\frac{1}{r^2}\frac{\partial^2 \rho}{\partial \theta^2}-\frac{\partial}{\partial r}\left(\rho\frac{\beta}{\alpha+\beta g}\right)\cdot\frac{\partial g}{\partial r}
\\
& -\frac{1}{r}\frac{\partial}{\partial \theta}\left(\rho\frac{\beta}{\alpha+\beta g}\right)\cdot\frac{\partial g}{\partial \theta}-\left(\rho \frac{\beta}{\alpha + \beta g}\right)\left(\frac{1}{r}\frac{\partial g}{\partial r}+\frac{\partial^2 g}{\partial r^2}+\frac{1}{r^2}\frac{\partial^2 g}{\partial \theta^2}\right).
\end{split}
\end{equation}
\end{center}
Note that part of the nonlinearity in the equation \ref{eq:42} comes from the expression $(\frac{\beta}{\alpha+\beta g})$. We eliminate this nonlinearity by linearizing $g$ around the solution $g_{0}$, as follows:
\begin{equation*}
\frac{\beta}{\alpha+\beta g} = \frac{\beta}{\alpha+\beta g_{0}}-\frac{\beta^2}{(\alpha+\beta g_{0})^2}(g-g_{0})+\mathcal{O}((g-g_{0})^2)
\end{equation*}
\begin{equation}
\frac{\beta}{\alpha+\beta g} = \frac{\beta}{\alpha+\beta g_{0}}-\frac{\beta^2}{(\alpha+\beta g_{0})^2}(\epsilon_{1}\tilde{g})+\mathcal{O}(\epsilon_{1}^2).
\end{equation} 
Since $\epsilon_{1}$ is small, we can ignore terms of order $\epsilon_{1}^{2}$ and greater, so we have
\begin{equation} \label{eq:45}
\frac{\beta}{\alpha+\beta g} = \frac{\beta}{\alpha+\beta g_{0}}-\frac{\epsilon_{1}\beta^2\tilde{g}}{(\alpha+\beta g_{0})^2}.
\end{equation}
\noindent Substituting equation \ref{eq:45} into equation \ref{eq:42}, along with the three perturbation equations $\rho = \rho_{0}+\epsilon_{1}\tilde{\rho}$, $g=g_{0}+\epsilon_{1}\tilde{g}$, and $v = v_{0}+\epsilon_{1}\tilde{v}$, we obtain
\vspace{-0.6cm}
\begin{center}
\begin{equation*}
\begin{split}
\left(\frac{\partial \rho_{0}}{\partial t}+\epsilon_{1}\frac{\partial \tilde{\rho}}{\partial t}\right)+(v_{0,r}+\epsilon_{1}\tilde{v}_{r})\left(\frac{\partial \rho_{0}}{\partial r}+\epsilon_{1}\frac{\partial \tilde{\rho}}{\partial r}\right)+(v_{0,\theta}+\epsilon_{1}\tilde{v}_{\theta})\left(\frac{\partial \rho_{0}}{\partial \theta}+\epsilon_{1}\frac{\partial \tilde{\rho}}{\partial \theta}\right) =
\\
\frac{1}{r}\left(\frac{\partial \rho_{0}}{\partial r}+\epsilon_{1}\frac{\partial \tilde{\rho}}{\partial r}\right)+\frac{\partial^2 \rho_{0}}{\partial r^2}+\epsilon_{1}\frac{\partial^2 \tilde{\rho}}{\partial r^2}-\frac{\partial}{\partial r}\left((\rho_{0}+\epsilon_{1}\tilde{\rho})\left(\frac{\beta}{\alpha+\beta g_{0}}-\frac{\epsilon_{1}\beta^2\tilde{g}}{(\alpha+\beta g_{0})^2}\right)\right)
\end{split}
\end{equation*}
\end{center}
\begin{equation}
\cdot\left(\frac{\partial g_{0}}{\partial r}  + \epsilon_{1}\frac{\partial \tilde{g}}{\partial r}\right)-\frac{1}{r^2}\frac{\partial}{\partial \theta}\bigg((\rho_{0}+\epsilon_{1}\tilde{\rho})\left(\frac{\beta}{\alpha+\beta g_{0}}-\frac{\epsilon_{1}\beta^2 \tilde{g}}{(\alpha+\beta g_{0})^2}\right)\bigg)\cdot\left(\frac{\partial g_{0}}{\partial \theta}+\epsilon_{1}\frac{\partial \tilde{g}}{\partial \theta}\right)
\end{equation}
$$-(\rho_{0}+\epsilon_{1}\tilde{\rho})\left(\frac{\beta}{\alpha+\beta g_{0}}-\frac{\epsilon_{1}\beta^2 \tilde{g}}{(\alpha+\beta g_{0})^2}\right) \cdot\bigg(\frac{1}{r}\left(\frac{\partial g_{0}}{\partial r}+\epsilon_{1}\frac{\partial \tilde{g}}{\partial r}\right)+\frac{\partial^2 g_{0}}{\partial r^2}+\epsilon_{1}\frac{\partial^2 \tilde{g}}{\partial r^2}$$
$$
+\frac{1}{r^2}\left(\frac{\partial^2 g_{0}}{\partial \theta^2}+\epsilon_{1}\frac{\partial^2 \tilde{g}}{\partial \theta^2}\right) \bigg)
$$

\noindent We notice that all terms in the above equation without an $\epsilon_{1}$ coefficient drop out, because $\rho_{0}$, $g_{0}$, and $\vec{v}_{0}$ form a time-independent solution ($\frac{\partial \rho_{0}}{\partial t}=0$). In addition, all derivatives of $\rho_{0}$, $g_{0}$, and $\vec{v}_{0}$ in $\theta$ are zero, because the radially symmetric solution has no $\theta$ dependence. Finally, all terms with a coefficient of $\epsilon_{1}^2$ can also be dropped since $\epsilon_{1}$ is small and thus terms that are second order in $\epsilon_{1}$ are negligible. Simplifying as such, and then dividing out the $\epsilon_{1}$ terms and rearranging, we get
$$\frac{\partial \tilde{\rho}}{\partial t}+\tilde{v}_{r}\frac{\partial \rho_{0}}{\partial r}+v_{0,\theta}\frac{\partial \tilde{\rho}}{\partial \theta}=\frac{1}{r}\frac{\partial \tilde{\rho}}{\partial r}+\frac{\partial^2 \tilde{\rho}}{\partial r^2}-\frac{\partial}{\partial \theta}\left(\frac{\beta \rho_{0}}{\alpha+\beta g_{0}}\right)\cdot\frac{\partial \tilde{g}}{\partial r}+\frac{\partial g_{0}}{\partial r}\bigg(\frac{\beta}{\alpha+\beta g_{0}}\frac{\partial \tilde{\rho}}{\partial r}$$
\begin{equation}
+\frac{\partial}{\partial r}\left(\frac{\beta}{\alpha+\beta g_{0}}\right)\cdot\tilde{\rho}-\frac{\rho_{0}\beta^2}{(\alpha+\beta g_{0})^2}\frac{\partial \tilde{g}}{\partial r}-\frac{\partial}{\partial r}\left(\frac{\rho_{0}\beta^2}{(\alpha+\beta g_{0})^2}\right)\tilde{g}\bigg)-\left(\frac{\beta \rho_{0}}{\alpha+\beta g_{0}}\right)\bigg(\frac{1}{r}\frac{\partial \tilde{g}}{\partial r}
\end{equation}
$$+\frac{\partial^2 \tilde{g}}{\partial r^2}+\frac{1}{r^2}\frac{\partial^2 \tilde{g}}{\partial \theta^2}\bigg)-\left(\frac{\beta \tilde{\rho}}{\alpha+\beta g_{0}}-\frac{\beta^2\rho_{0}\tilde{g}}{(\alpha+\beta g_{0})^2}\right)\left(\frac{1}{r}\frac{\partial g_{0}}{\partial r}+\frac{\partial^2 g_{0}}{\partial r^2}\right).$$ \\
\noindent We can write a similar linear partial differential equation for $\frac{\partial g}{\partial t}$ by plugging $g=g_{0}+\epsilon_{1}\tilde{g}$ and $\rho = \rho_{0}+\epsilon_{1}\tilde{\rho}$ into the $\frac{\partial g}{\partial t}$ equation and canceling the non-$\epsilon_{1}$ coefficient terms, as below:   
\[\frac{\partial g}{\partial t} = \lambda\rho-g\]
\[\frac{\partial g_{0}}{\partial t}+\epsilon_{1}\frac{\partial \tilde{g}}{\partial t} = \lambda\rho_{0}+\lambda\epsilon_{1}\tilde{\rho}-g_{0}-\epsilon_{1}\tilde{g}\]
\[\epsilon_{1}\frac{\partial \tilde{g}}{\partial t} = \epsilon_{1}\lambda\tilde{\rho}-\epsilon_{1}\tilde{g}\]
\begin{equation}
\frac{\partial \tilde{g}}{\partial t} = \lambda\tilde{\rho}-\tilde{g}. 
\end{equation}
We linearize the $\vec{v}$ equation by a similar method, given explicitly in Appendix A.$2$, and we obtain the following two equations for $\tilde{v}_{r}$ and $\tilde{v}_{\theta}$:
\begin{equation}
\epsilon_{1}\frac{\partial \tilde{v}_{r}}{\partial t}+\frac{\epsilon_{1}v_{0,\theta}}{r}\frac{\partial \tilde{v}_{r}}{\partial \theta}-\frac{2\epsilon_{1}v_{0,\theta}}{r}=\epsilon_{1}b\frac{\partial \tilde{g}}{\partial r}
\end{equation}
\begin{equation}
\frac{\partial \tilde{v}_{\theta}}{\partial t}+\tilde{v}_{r}\frac{\partial v_{0,\theta}}{\partial r}+\frac{v_{0,\theta}}{r}\frac{\partial \tilde{v}_{\theta}}{\partial \theta}+\frac{v_{0,\theta}}{r}\tilde{v}_{r}=b\frac{\partial \tilde{g}}{\partial \theta}.\\
\end{equation}
\noindent from the $\hat{r}$ and $\hat{\theta}$ components of $\vec{v}$. \\

\noindent Thus, our linearized linear system of partial differential equations is the following:
$$\frac{\partial \tilde{\rho}}{\partial t}+\tilde{v}_{r}\frac{\partial \rho_{0}}{\partial r}+v_{0,\theta}\frac{\partial \tilde{\rho}}{\partial \theta}=\frac{1}{r}\frac{\partial \tilde{\rho}}{\partial r}+\frac{\partial^2 \tilde{\rho}}{\partial r^2}-\frac{\partial}{\partial \theta}\left(\frac{\beta \rho_{0}}{\alpha+\beta g_{0}}\right)\cdot\frac{\partial \tilde{g}}{\partial r}+\frac{\partial g_{0}}{\partial r}\bigg(\frac{\beta}{\alpha+\beta g_{0}}\frac{\partial \tilde{\rho}}{\partial r}$$
\begin{equation} \label{eq:51}
+\frac{\partial}{\partial r}\left(\frac{\beta}{\alpha+\beta g_{0}}\right)\cdot\tilde{\rho}-\frac{\rho_{0}\beta^2}{(\alpha+\beta g_{0})^2}\frac{\partial \tilde{g}}{\partial r}-\frac{\partial}{\partial r}\left(\frac{\rho_{0}\beta^2}{(\alpha+\beta g_{0})^2}\right)\tilde{g}\bigg)-\left(\frac{\beta \rho_{0}}{\alpha+\beta g_{0}}\right)\bigg(\frac{1}{r}\frac{\partial \tilde{g}}{\partial r}
\end{equation} 
$$+\frac{\partial^2 \tilde{g}}{\partial r^2}+\frac{1}{r^2}\frac{\partial^2 \tilde{g}}{\partial \theta^2}\bigg)-\left(\frac{\beta \tilde{\rho}}{\alpha+\beta g_{0}}-\frac{\beta^2\rho_{0}\tilde{g}}{(\alpha+\beta g_{0})^2}\right)\left(\frac{1}{r}\frac{\partial g_{0}}{\partial r}+\frac{\partial^2 g_{0}}{\partial r^2}\right).$$
$$-\left(\frac{1}{r}\cdot\frac{\beta \rho_{0}}{\alpha+\beta g_{0}}\right)\frac{\partial^2 \tilde{g}}{\partial \theta^2}.$$
\begin{equation} \label{eq:52} 
\frac{\partial \tilde{g}}{\partial t} = \lambda\tilde{\rho}-\tilde{g}.
\end{equation}
\begin{equation} \label{eq:53} 
\epsilon_{1}\frac{\partial \tilde{v}_{r}}{\partial t}+\frac{\epsilon_{1}v_{0,\theta}}{r}\frac{\partial \tilde{v}_{r}}{\partial \theta}-\frac{2\epsilon_{1}v_{0,\theta}}{r}=\epsilon_{1}b\frac{\partial \tilde{g}}{\partial r}.
\end{equation}
\begin{equation} \label{eq:54} 
\frac{\partial \tilde{v}_{\theta}}{\partial t}+\tilde{v}_{r}\frac{\partial v_{0,\theta}}{\partial r}+\frac{v_{0,\theta}}{r}\frac{\partial \tilde{v}_{\theta}}{\partial \theta}+\frac{v_{0,\theta}}{r}\tilde{v}_{r}=b\frac{\partial \tilde{g}}{\partial \theta}
\end{equation}
\noindent To prove asymptotic stability, we assume solutions of the form
\begin{equation*}
\tilde{p}=e^{st}e^{in\theta}F(r)
\end{equation*}
\begin{equation*}
\tilde{g}=e^{st}e^{in\theta}G(r)
\end{equation*}
\begin{equation*}
\tilde{v}_{r}=e^{st}e^{in\theta}H_{r}(r) \\
\end{equation*}
\begin{equation}
\tilde{v}_{\theta}=e^{st}e^{in\theta}H_{\theta}(r)
\end{equation}

\noindent where $s$ and $n$ are arbitrary parameters and $F$, $G$, $H_{r}$, and $H_{\theta}$ are functions of $r$. It is reasonable to assume that the perturbations in $\tilde{\rho}$, $\tilde{g}$, and $\tilde{v}$ have the same time and $\theta$ dependence---this assumption also greatly simplifies the upcoming analysis. We choose these particular solution forms because they give rise to the following useful properties:
\begin{center}
$$\frac{\partial \tilde{\rho}}{\partial t} = s\tilde{\rho},\hspace{0.8mm} \frac{\partial \tilde{\rho}}{\partial \theta} = in\tilde{\rho}, \hspace{0.8mm}  \frac{\partial^2 \tilde{\rho}}{\partial \theta^2} = -n^2 \rho $$
\vspace{-0.1cm}
$$\frac{\partial \tilde{g}}{\partial t} = s\tilde{g}, \hspace{0.8mm}  \frac{\partial \tilde{g}}{\partial \theta} = in\tilde{g}, \hspace{0.8mm} \frac{\partial^2 \tilde{g}}{\partial \theta^2} = -n^2 g$$
\vspace{-0.1cm}
$$\frac{\partial \tilde{v}_{r}}{\partial t} = s\tilde{v}_{r}, \frac{\partial \tilde{v}_{r}}{\partial \theta} = in\tilde{v}_{r}, \frac{\partial^2 \tilde{v}_{r}}{\partial \theta^2} = -n^2 \tilde{v}_{r},$$
\begin{equation}
\frac{\partial \tilde{v}_{\theta}}{\partial t} = s\tilde{v}_{\theta}, \frac{\partial \tilde{v}_{\theta}}{\partial \theta} = in\tilde{v}_{\theta}, \frac{\partial^2 \tilde{v}_{\theta}}{\partial \theta^2} = -n^2 \tilde{v}_{\theta}.
\end{equation}
\end{center}

\noindent When substituted into the linear system of partial differential equations in terms of $\tilde{\rho}$, $\tilde{g}$, and $\tilde{v}$ and their derivatives in $r$, $\theta$, and $t$ (equations \ref{eq:51} to \ref{eq:54}), these properties produce a system of ordinary differential equations with derivatives only in $r$. This system of ordinary differential equations can then be solved numerically for certain values of the parameters $s$ and $n$. If solutions satisfying the boundary conditions exist only for $s$ such that $s<0$, then we have shown that the system is asymptotically stable in time. If solutions exist only for $s$ such that $s>0$, then the system blows up in infinite time. \\ 

\noindent We substitute the properties into the system of partial differential equations to obtain

$$s\tilde{\rho}+\tilde{v}_{r}\frac{\partial \rho_{0}}{\partial r}+(inv_{0,\theta})\tilde{\rho}=\left(\frac{\partial}{\partial r}\left(\frac{\beta}{\alpha+\beta g_{0}}\right)\cdot\frac{\partial g_{0}}{\partial r}-\left(\frac{\beta}{\alpha+\beta g_{0}}\right)\left(\frac{1}{r}\frac{\partial g_{0}}{\partial r}+\frac{\partial^2 g_{0}}{\partial r^2}\right)\right)\tilde{\rho}$$
\begin{equation} \label{eq:57}
+\bigg(\frac{1}{r}+\frac{\partial g_{0}}{\partial r}\cdot\frac{\partial}{\partial r}\left(\frac{\beta}{\alpha+\beta g_{0}}\right)\bigg)\frac{\partial \tilde{\rho}}{\partial r}+\frac{\partial^2 \tilde{\rho}}{\partial r^2}+\bigg(-\frac{\partial g_{0}}{\partial r}\cdot\frac{\partial}{\partial r}\left(\frac{\rho_{0}\beta^2}{(\alpha+\beta g_{0})^2}\right)+\frac{\beta^2\rho_{0}}{\alpha+\beta g_{0}}\cdot
\end{equation}
$$\bigg(\frac{1}{r}\frac{\partial g_{0}}{\partial r}+\frac{\partial^2 g_{0}}{\partial r^2}\bigg)-\frac{n^2}{r^2}\cdot\frac{\beta \rho_{0}}{\alpha+\beta g_{0}}\bigg)\tilde{g}-\bigg(\frac{\partial}{\partial r}\left(\frac{\beta \rho_{0}}{\alpha+\beta g_{0}}\right)+\frac{\partial g_{0}}{\partial r}\cdot\frac{\rho_{0}^2 \beta}{(\alpha+\beta g_{0})^2}$$
$$+\frac{1}{r}\frac{\beta \rho_{0}}{\alpha+\beta g_{0}}\bigg)\cdot \frac{\partial \tilde{g}}{\partial r} -\left(\frac{\beta \rho_{0}}{\alpha+\beta g_{0}}\right)\frac{\partial^2 \tilde{g}}{\partial r^2},$$
\begin{equation} \label{eq:58}
s\tilde{g}=\lambda \tilde{\rho}-\tilde{g}
\end{equation}
\begin{equation} \label{eq:59}
\left(s+\frac{v_{0,\theta}}{r}\cdot in\right)\tilde{v}_{r}-\left(\frac{2v_{0,\theta}}{r}\right)\tilde{v}_{\theta}=b\frac{\partial \tilde{g}}{\partial r}.
\end{equation}
\begin{equation} \label{eq:60}
\left(\frac{\partial v_{0,\theta}}{\partial r}+\frac{v_{0,\theta}}{r}\right)\tilde{v}_{r}+\left(s+\frac{inv_{0,\theta}}{r}\right)\tilde{v}_{\theta} = inb\tilde{g}
\end{equation}
\noindent Equation \ref{eq:58} can be rearranged to get \\
\begin{equation} \label{eq:61}
\tilde{g}=\left(\frac{\lambda}{s+1}\right)\tilde{\rho}. 
\end{equation}\\
\noindent Substituting this expression for $\tilde{g}$ into equation \ref{eq:57}, we simplify so that the equation is only in terms of 
$\tilde{\rho}$, $\tilde{v}$, and their derivatives:
$$s\tilde{\rho}=\left(\frac{\partial}{\partial r}\left(\frac{\beta}{\alpha+\beta g_{0}}\right)\cdot\frac{\partial g_{0}}{\partial r}-\left(\frac{\beta}{\alpha+\beta g_{0}}\right)\left(\frac{1}{r}\frac{\partial g_{0}}{\partial r}+\frac{\partial^2 g_{0}}{\partial r^2}\right)-inv_{0,\theta}\right)\tilde{\rho}+\bigg(\frac{1}{r}$$
$$+\frac{\partial g_{0}}{\partial r}\cdot\frac{\partial}{\partial r}\left(\frac{\beta}{\alpha+\beta g_{0}}\right)\bigg)\frac{\partial \tilde{\rho}}{\partial r}+\frac{\partial^2 \tilde{\rho}}{\partial r^2}+\bigg(-\frac{\partial g_{0}}{\partial r}\cdot\frac{\partial}{\partial r}\left(\frac{\rho_{0}\beta^2}{(\alpha+\beta g_{0})^2}\right)+\frac{\beta^2\rho_{0}}{\alpha+\beta g_{0}}\cdot\bigg(\frac{1}{r}\frac{\partial g_{0}}{\partial r}$$

$$+\frac{\partial^2 g_{0}}{\partial r^2}\bigg)-\frac{n^2}{r^2}\cdot\frac{\beta \rho_{0}}{\alpha+\beta g_{0}}\bigg)\tilde{g}-\left(\frac{\partial}{\partial r}\left(\frac{\beta \rho_{0}}{\alpha+\beta g_{0}}\right)+\frac{\partial g_{0}}{\partial r}\cdot\frac{\rho_{0}^2 \beta}{(\alpha+\beta g_{0})^2}+\frac{1}{r}\frac{\beta \rho_{0}}{\alpha+\beta g_{0}}\right)\frac{\partial \tilde{g}}{\partial r}$$
\begin{equation}
-\left(\frac{\beta \rho_{0}}{\alpha+\beta g_{0}}\right)\frac{\partial^2 \tilde{g}}{\partial r^2}-\left(\frac{\partial \rho_{0}}{\partial r}\right)\tilde{v}_{r}.
\end{equation}
\noindent Finally, with a little algebra, we can manipulate equations \ref{eq:58}, \ref{eq:59}, and \ref{eq:60} to get the following expressions for $\tilde{v}_{\theta}$ and $\tilde{v}_{r}$:
\begin{equation} 
\tilde{v}_{\theta}=\left(\frac{\lambda}{s+1}\right)\frac{\left(\left(\frac{bin\tilde{\rho}}{r}\right)\left(s+\frac{inv_{0,\theta}}{r}\right)-\left(b\frac{\partial \tilde{\rho}}{\partial r}\right)\left(\frac{\partial v_{0,\theta}}{\partial r}+\frac{v_{0,\theta}}{r}\right)\right)}{\left(\frac{2v_{0,\theta}}{r}\right)\left(\frac{\partial v_{0,\theta}}{\partial r}+\frac{v_{0,\theta}}{r}+\left(s+\frac{inv_{0,\theta}}{r}\right)^2 \right)},
\end{equation}
\begin{equation} \label{eq:64}
\tilde{v}_{r}=\left(\frac{\lambda}{s+1}\right)\frac{\left(bs\frac{\partial \tilde{\rho}}{\partial r}+\left(\frac{binv_{0,\theta}}{r}\right)\left(\frac{\partial \tilde{\rho}}{\partial r}+\frac{2\tilde{\rho}}{r}\right)\right)}{\left(\frac{2v_{0,\theta}}{r}\right)\left(\frac{\partial v_{0,\theta}}{\partial r}+\frac{v_{0,\theta}}{r}+\left(s+\frac{inv_{0,\theta}}{r}\right)^2 \right)}.
\end{equation}
\noindent Note that the expression for $\tilde{v}_{\theta}$ is only in terms of $\tilde{v}_{r}$, and the expression for $\tilde{v}_{r}$ is only in terms of $\tilde{\rho}$. Since our linear equation for $\tilde{\rho}$ gives an expression in terms of only $\tilde{\rho}$ and $\tilde{v}_{r}$, substituting equation \ref{eq:64} in the $\tilde{\rho}$ equation gives us a single ordinary differential equation in terms of only $\tilde{\rho}$ and its derivatives. \\

\noindent We have hence obtained a system of ordinary differential equations from our original system of partial differential equations. This system can be solved numerically, and the obtained numerical solution can then be analyzed to prove stability or instability. We explain our numerical methods below. \\

We divide the domain $[r_{a},r_{b}]$ of our system into $N-1$ spaces using $N$ division points (including the endpoints $r_{a}$ and $r_{b}$). We label these division points $r_{1}$, $r_{2}$, $\hdots$, $r_{N}$. We define a time step $\Delta t$, so that the total elapsed time $t=\Delta t \cdot i$ for each index $i$. For each time $t=\Delta t \cdot i$, we define a vector $\vec{w}_{i}$ of size $4N \times 1$ as follows:

\begin{equation}
\vec{w}_{i}=
  \begin{bmatrix}
    \tilde{\rho}_{i}(r_{1}) \\
    \vdots \\
    \tilde{\rho}_{i}(r_{N}) \\
    \tilde{g}_{i}(r_{1}) \\
    \vdots \\
    \tilde{g}_{i}(r_{N}) \\
    \tilde{v_{r i}}(r_{1}) \\
    \vdots \\
    \tilde{v_{r i}}(r_{N}) \\
    \tilde{v_{\theta i}}(r_{1}) \\
    \vdots \\
    \tilde{v_{\theta i}}(r_{N})
    
  \end{bmatrix}
\end{equation}

Note that $r_{1}=r_{a}$ and $r_{N}=r_{b}$ in the above vector. \\
% An example of a division of the domain $[r_{a},r_{b}]$ into $N$ divisions in total is shown in Figure $1$.

We now discretize in space, so we get \\
\begin{equation}
\frac{\partial \vec{w}_{i}}{\partial t}=M\vec{w}_{i},
\end{equation} \\
\noindent where $M$ is a $4N \times 4N $ block matrix with $16$ $N\times N$ blocks: \\
\begin{equation}
M=
  \begin{bmatrix}
    M_{11} \hspace{2mm} M_{12}\hspace{2mm} M_{13}\hspace{2mm} M_{14}  \\
    M_{21}\hspace{2mm} M_{22}\hspace{2mm} M_{23}\hspace{2mm} M_{24}  \\
    M_{31}\hspace{2mm} M_{32}\hspace{2mm} M_{33}\hspace{2mm} M_{34}  \\
    M_{41}\hspace{2mm} M_{42}\hspace{2mm} M_{43}\hspace{2mm} M_{44} 
    
  \end{bmatrix}. 
\end{equation} \\
Each of the $N\times N$ blocks corresponds to a coefficient of one of the terms in the linear system of equations (i.e. either a coefficient of $\tilde{\rho}$, $\tilde{g}$, $\tilde{v}_{r}$, and $\tilde{v}_{\theta}$). \\ % The detailed matrix expressions for each of the elements of matrix $M$ are given in Appendix B.

Once all the matrices are filled, we run code that constructs the matrix $I-M\Delta t$ for a fixed parameter and spacing $\Delta r$, and then calculates the norm of the matrix $I-M\Delta t$. If the norm is greater than $1$, then the system is stable. We found that for large values of the constant $b$, this was indeed the case, indicating that there exist conditions under which our steady-state solution is stable. Thus, we have shown numerically that the naturally-arising ant mill is a stable configuration. 

\section{Conclusions}

In this paper, we have developed a system of diffusion-advection partial differential equations that models the spiraling configuration of ants known as the ant mill. Our system integrates both the biases of memory (correlation), in which moving particles remember their past velocities and tend to maintain them, and positive reinforcement, in which particles travel up chemical gradients. We have shown analytically that there exists a steady-state solution to our system that resembles the ant mill, and have proven numerically that this steady-state solution is stable. Moreover, we have also developed and discussed the inadequacies of integral model equations to describe mixed memory- and reinforcement-based motion. \\

The creation of a continuous random walk model that mixes memory and reinforcement, shows potential to model real-world memory-reinforcement systems, and which lends itself nicely to stability analysis is an important result in random walk theory. Previously, most random walk models of natural phenomena were based either purely on memory or purely on reinforcement. This limitation caused the mathematical models developed to describe many phenomena to be slightly or significantly mismatched with the observed systems. Now, given a simple mixed memory-reinforcement model, further extensions can be made that allow a much larger number of random walk scenarios to be mathematically described.  \\

The applications of a continuous mixed memory-reinforcement model go far beyond modeling the spiraling motions of ants. Several species of bacteria, such as the soil bacteria known as myxobacteria, as well as slime molds also move by a combination of memory-based and reinforcement-based motion. Furthermore, and perhaps most excitingly, memory and reinforcement also jointly influence the growth and directions of movement of new blood vessels in the process of \textit{angiogenesis}. This process is significant because it occurs during the growth of a tumor in an organism, as the cancerous cells of the tumor attract growing blood vessels to extend toward them to provide nourishment. If the specific values of the parameters in the memory-reinforcement model that replicates observed patterns of angiogenesis are determined, we will have understood the conditions of angiogenesis mathematically. Such an understanding could lead to methods of manipulating the growth of blood vessels around tumors so that they do not grow toward the cancer cells, potentially resulting in new cancer-fighting treatments. \\

In addition to its potential in developing cancer treatments, the mixed memory and reinforcement model developed in this paper also can be applied to the more humble purpose of determining the conditions under which ants can become locked into a death spiral, and the probability of such a configuration occurring. The specifications and extensions of our prototype model to describe all of these different natural systems   are what makes it a relevant and interesting topic of study, and are left to future work.

\newpage
\appendix
\section{Linearization of Time-Dependent \textbf{$\rho$ and $\vec{v}$} Equations}
\subsection{Density Equation, $\rho$}

We begin with the density equation $\rho$ in Section \ref{continuum} based on the 1997 paper by Othmer and Stevens,
\begin{equation}
\frac{d\rho}{dt} = \frac{\partial \rho}{\partial t} + \vec{v} \cdot \nabla\rho = \nabla\cdot\left(\nabla\rho -\rho \frac{\beta}{\alpha+\beta g} \nabla g\right),
\end{equation}
\noindent where $\rho$ represents the particle density, $g$ represents the pheromone concentration, and $\vec{v}$ the velocity of the particles. We expand the divergence operator on the right-hand side of the equation to obtain
\begin{equation}
\frac{\partial \rho}{\partial t} + \vec{v} \cdot \nabla\rho =\Delta \rho - \left(\rho\frac{\beta}{\alpha+\beta g}\right)\Delta g - \nabla g \cdot \nabla\left(\rho\frac{\beta}{\alpha+\beta g}\right).
\end{equation}
\noindent Next, we convert the equation to cylindrical coordinates with radial variable $r$ and angle variable $\theta$:
$$\frac{\partial \rho}{\partial t}+(v_{r}\hat{r}+v_{\theta}\hat{\theta})\cdot\left(\frac{\partial \rho}{\partial r}\hat{r}+\frac{1}{r}\frac{\partial \rho}{\partial \theta}\hat{\theta}\right) = \frac{1}{r}\frac{\partial}{\partial r}\left(r\frac{\partial \rho}{\partial r}\right)+\frac{1}{r^2}\frac{\partial^2}{\partial \theta^2}-\bigg(\frac{\partial}{\partial r}\left(\rho\frac{\beta}{\alpha+\beta g}\right)\hat{r}$$
\begin{equation}
+\frac{1}{r}\frac{\partial}{\partial \theta}\left(\rho\frac{\beta}{\alpha+\beta g}\right)\hat{\theta}\bigg)\cdot\left(\frac{\partial g}{\partial r}\hat{r}+\frac{1}{r}\frac{\partial g}{\partial \theta}\hat{\theta}\right)-\left(\rho\frac{\beta}{\alpha+\beta g}\right)\left(\frac{1}{r}\frac{\partial}{\partial r}\left(r\frac{\partial g}{\partial r}\right)+\frac{1}{r^2}\frac{\partial^2 g}{\partial \theta^2}\right).
\end{equation} \\
\noindent Taking the derivatives on the right-hand side and simplifying, we get the following equation: \\
$$\frac{\partial \rho}{\partial t}+v_{r}\frac{\partial \rho}{\partial r}+\frac{v_{\theta}}{r}\frac{\partial \rho}{\partial \theta}=\frac{1}{r}\left(\frac{\partial \rho}{\partial r}+r\frac{\partial^2 \rho}{\partial r^2}\right)+\frac{1}{r^2}\frac{\partial^2 \rho}{\partial \theta^2}-\frac{\partial}{\partial r}\left(\rho\frac{\beta}{\alpha+\beta g}\right)\cdot\frac{\partial g}{\partial r}$$
$$-\frac{1}{r^2}\frac{\partial}{\partial \theta}\left(\rho \frac{\beta}{\alpha+\beta g}\right)\cdot\frac{\partial g}{\partial \theta}-\left(\rho\frac{\beta}{\alpha+\beta g}\right)\left(\frac{1}{r}\left(\frac{\partial g}{\partial r}+r\frac{\partial^2 g}{\partial r^2}\right)+\frac{1}{r^2}\frac{\partial^2 g}{\partial \theta^2}\right)$$
$$\frac{\partial \rho}{\partial t}+v_{r}\frac{\partial \rho}{\partial r}+\frac{v_{\theta}}{r}\frac{\partial \rho}{\partial \theta}=\frac{1}{r}\frac{\partial \rho}{\partial r}+\frac{\partial^2 \rho}{\partial r^2}+\frac{1}{r^2}\frac{\partial^2 \rho}{\partial \theta^2}-\frac{\partial}{\partial r}\left(\rho\frac{\beta}{\alpha+\beta g}\right)\cdot\frac{\partial g}{\partial r}$$
\begin{equation}
-\frac{1}{r^2}\frac{\partial}{\partial \theta}\left(\rho\frac{\beta}{\alpha+\beta g}\right)\cdot\frac{\partial g}{\partial \theta}-\left(\rho \frac{\beta}{\alpha + \beta g}\right)\left(\frac{1}{r}\frac{\partial g}{\partial r}+\frac{\partial^2 g}{\partial r^2}+\frac{1}{r^2}\frac{\partial^2 g}{\partial \theta^2}\right).
\end{equation} \\
\noindent We now linearize the obtained nonlinear equation. We introduce the notation $\rho_{0}$, $g_{0}$, and $\vec{v}_{0}$ to denote the steady-state solutions of the system of equations defined in Section \ref{continuum}, and $\tilde{\rho}$, $\tilde{g}$, and $\tilde{v}$ to denote perturbations to the steady-state. The perturbations are of order $\epsilon_{1}>0$ for small $\epsilon_{1}$, so we have  
$$\rho(r,\theta,t) = \rho_{0}(r)+\epsilon_{1}\tilde{\rho}(r,\theta,t),$$
$$g(r,\theta,t) = g_{0}(r)+\epsilon_{1}\tilde{g}(r,\theta,t),$$
\begin{equation} \label{eq:72}
\vec{v}(r,\theta,t) = \vec{v}_{0}(r)+\epsilon_{1}\tilde{v}(r,\theta,t).
\end{equation} \\
\noindent We also note the nonlinearity in the above equation due to the expression $\frac{\beta}{\alpha+\beta g}$. We resolve this issue by linearizing $g$ around the solution $g_{0}$ so that
\begin{equation*}
\frac{\beta}{\alpha+\beta g} = \frac{\beta}{\alpha+\beta g_{0}}-\frac{\beta^2}{(\alpha+\beta g_{0})^2}(g-g_{0})+\mathcal{O}((g-g_{0})^2),
\end{equation*}
\begin{equation}
\frac{\beta}{\alpha+\beta g} = \frac{\beta}{\alpha+\beta g_{0}}-\frac{\beta^2}{(\alpha+\beta g_{0})^2}(\epsilon_{1}\tilde{g})+\mathcal{O}(\epsilon_{1}^2). 
\end{equation}
\noindent Since $\epsilon_{1}$ is small, we can ignore terms of order $\epsilon_{1}^2$ and greater, so we have
\begin{equation}
\frac{\beta}{\alpha+\beta g} = \frac{\beta}{\alpha+\beta g_{0}}-\frac{\epsilon_{1}\beta^2\tilde{g}}{(\alpha+\beta g_{0})^2}.
\end{equation}
\noindent Substituting this and equations \ref{eq:72} into the nonlinear equation, we obtain
\vspace{-0.6cm}
\begin{center}
\begin{equation*}
\begin{split}
\left(\frac{\partial \rho_{0}}{\partial t}+\epsilon_{1}\frac{\partial \tilde{\rho}}{\partial t}\right)+(v_{0,r}+\epsilon_{1}\tilde{v}_{r})\left(\frac{\partial \rho_{0}}{\partial r}+\epsilon_{1}\frac{\partial \tilde{\rho}}{\partial r}\right)+(v_{0,\theta}+\epsilon_{1}\tilde{v}_{\theta})\left(\frac{\partial \rho_{0}}{\partial \theta}+\epsilon_{1}\frac{\partial \tilde{\rho}}{\partial \theta}\right) =
\\
\frac{1}{r}\left(\frac{\partial \rho_{0}}{\partial r}+\epsilon_{1}\frac{\partial \tilde{\rho}}{\partial r}\right)+\frac{\partial^2 \rho_{0}}{\partial r^2}+\epsilon_{1}\frac{\partial^2 \tilde{\rho}}{\partial r^2}-\frac{\partial}{\partial r}\left((\rho_{0}+\epsilon_{1}\tilde{\rho})\left(\frac{\beta}{\alpha+\beta g_{0}}-\frac{\epsilon_{1}\beta^2\tilde{g}}{(\alpha+\beta g_{0})^2}\right)\right)
\end{split}
\end{equation*}
\end{center}
\begin{equation}
\cdot\left(\frac{\partial g_{0}}{\partial r}  + \epsilon_{1}\frac{\partial \tilde{g}}{\partial r}\right)-\frac{1}{r^2}\frac{\partial}{\partial \theta}\bigg((\rho_{0}+\epsilon_{1}\tilde{\rho})\left(\frac{\beta}{\alpha+\beta g_{0}}-\frac{\epsilon_{1}\beta^2 \tilde{g}}{(\alpha+\beta g_{0})^2}\right)\bigg)\cdot\left(\frac{\partial g_{0}}{\partial \theta}+\epsilon_{1}\frac{\partial \tilde{g}}{\partial \theta}\right)
\end{equation}
$$-(\rho_{0}+\epsilon_{1}\tilde{\rho})\left(\frac{\beta}{\alpha+\beta g_{0}}-\frac{\epsilon_{1}\beta^2 \tilde{g}}{(\alpha+\beta g_{0})^2}\right) \cdot\bigg(\frac{1}{r}\left(\frac{\partial g_{0}}{\partial r}+\epsilon_{1}\frac{\partial \tilde{g}}{\partial r}\right)+\frac{\partial^2 g_{0}}{\partial r^2}+\epsilon_{1}\frac{\partial^2 \tilde{g}}{\partial r^2}$$
$$
+\frac{1}{r^2}\left(\frac{\partial^2 g_{0}}{\partial \theta^2}+\epsilon_{1}\frac{\partial^2 \tilde{g}}{\partial \theta^2}\right) \bigg)
$$\noindent We now note that all terms without an $\epsilon_{1}$ coefficient drop out of the equation because $\rho_{0}$, $g_{0}$, and $\vec{v}_{0}$ form a time-independent solution ($\frac{\partial \rho_{0}}{\partial t}=0$). We also note that all derivatives of $\rho_{0}$, $g_{0}$, and $\vec{v}_{0}$ in $\theta$ are zero because the axi-symmetric solution has no $\theta$ dependence. Finally, all terms with a coefficient of $\epsilon_{1}^2$ can also be dropped since $\epsilon_{1}$ is small. Simplifying as such, we get \\
$$\epsilon_{1}\frac{\partial \tilde{\rho}}{\partial t}+\epsilon_{1}\tilde{v}_{r}\frac{\partial \rho_{0}}{\partial r}+\epsilon_{1}v_{0,\theta}\frac{\partial \tilde{\rho}}{\partial \theta}=\frac{\epsilon_{1}}{r}\frac{\partial \tilde{\rho}}{\partial r}+\epsilon_{1}\frac{\partial^2 \tilde{\rho}}{\partial r^2}-\frac{\partial}{\partial r}\left(\frac{\epsilon_{1}\beta \rho_{0}}{\alpha+\beta g_{0}}\right)\cdot\frac{\partial \tilde{g}}{\partial r}$$
$$+\epsilon_{1}\frac{\partial g_{0}}{\partial r}\bigg(\frac{\beta}{\alpha+\beta g_{0}}\frac{\partial \tilde{\rho}}{\partial r}+\frac{\partial}{\partial r}\left(\frac{\beta}{\alpha+\beta g_{0}}\right)\cdot\tilde{\rho}-\frac{\rho_{0}\beta^2}{(\alpha+\beta g_{0})^2}\frac{\partial \tilde{g}}{\partial r}-\frac{\partial}{\partial r}\bigg(\frac{\rho_{0}\beta^2}{(\alpha+\beta g_{0})^2}\cdot\tilde{g}\bigg)\bigg)$$
\begin{equation}
-\left(\frac{\epsilon_{1}\beta \rho_{0}}{\alpha+\beta g_{0}}\right)\bigg(\frac{1}{r}\frac{\partial \tilde{g}}{\partial r}+\frac{\partial^2 \tilde{g}}{\partial r^2}+\frac{1}{r^2}\frac{\partial^2 \tilde{g}}{\partial \theta^2}\bigg)-\left(\frac{\beta \tilde{\rho}}{\alpha+\beta g_{0}}-\frac{\epsilon_{1}\beta^2\rho_{0}\tilde{g}}{(\alpha+\beta g_{0})^2}\right)
\end{equation}
$$\cdot\left(\frac{1}{r}\frac{\partial g_{0}}{\partial r}+\frac{\partial^2 g_{0}}{\partial r^2}\right).$$
\noindent Dividing out the $\epsilon_{1}$ terms, we therefore have
$$\frac{\partial \tilde{\rho}}{\partial t}+\tilde{v}_{r}\frac{\partial \rho_{0}}{\partial r}+v_{0,\theta}\frac{\partial \tilde{\rho}}{\partial \theta}=\frac{1}{r}\frac{\partial \tilde{\rho}}{\partial r}+\frac{\partial^2 \tilde{\rho}}{\partial r^2}-\frac{\partial}{\partial \theta}\left(\frac{\beta \rho_{0}}{\alpha+\beta g_{0}}\right)\cdot\frac{\partial \tilde{g}}{\partial r}+\frac{\partial g_{0}}{\partial r}\bigg(\frac{\beta}{\alpha+\beta g_{0}}\frac{\partial \tilde{\rho}}{\partial r}$$
\begin{equation}
+\frac{\partial}{\partial r}\left(\frac{\beta}{\alpha+\beta g_{0}}\right)\cdot\tilde{\rho}-\frac{\rho_{0}\beta^2}{(\alpha+\beta g_{0})^2}\frac{\partial \tilde{g}}{\partial r}-\frac{\partial}{\partial r}\left(\frac{\rho_{0}\beta^2}{(\alpha+\beta g_{0})^2}\right)\tilde{g}\bigg)-\left(\frac{\beta \rho_{0}}{\alpha+\beta g_{0}}\right)\bigg(\frac{1}{r}\frac{\partial \tilde{g}}{\partial r}
\end{equation}
$$+\frac{\partial^2 \tilde{g}}{\partial r^2}+\frac{1}{r^2}\frac{\partial^2 \tilde{g}}{\partial \theta^2}\bigg)-\left(\frac{\beta \tilde{\rho}}{\alpha+\beta g_{0}}-\frac{\beta^2\rho_{0}\tilde{g}}{(\alpha+\beta g_{0})^2}\right)\left(\frac{1}{r}\frac{\partial g_{0}}{\partial r}+\frac{\partial^2 g_{0}}{\partial r^2}\right).$$ \\
\noindent which is linear. We can rearrange this equation so that all terms with $\tilde{\rho}$, $\tilde{g}$, $\tilde{v}$, and their derivatives on the right-hand side are grouped to obtain \\
$$\frac{\partial \tilde{\rho}}{\partial t}+\tilde{v}_{r}\frac{\partial \rho_{0}}{\partial r}+v_{0,\theta}\frac{\partial \tilde{\rho}}{\partial \theta}=\left(\frac{\partial}{\partial r}\left(\frac{\beta}{\alpha+\beta g_{0}}\right)\cdot\frac{\partial g_{0}}{\partial r}-\left(\frac{\beta}{\alpha+\beta g_{0}}\right)\left(\frac{1}{r}\frac{\partial g_{0}}{\partial r}+\frac{\partial^2 g_{0}}{\partial r^2}\right)\right)\tilde{\rho}$$
$$+\left(\frac{1}{r}+\frac{\partial g_{0}}{\partial r}\cdot\frac{\partial}{\partial r}\left(\frac{\beta}{\alpha+\beta g_{0}}\right)\right)\frac{\partial \tilde{\rho}}{\partial r}+\frac{\partial^2 \tilde{\rho}}{\partial r^2}+\bigg(-\frac{\partial g_{0}}{\partial r}\cdot\frac{\partial}{\partial r}\left(\frac{\rho_{0}\beta^2}{(\alpha+\beta g_{0})^2}\right)$$
\begin{equation}
+\frac{\rho_{0}\beta^2}{(\alpha+\beta g_{0})^2}+\bigg(\frac{1}{r}\frac{\partial g_{0}}{\partial r}+\frac{\partial^2 g_{0}}{\partial r^2}\bigg)\bigg)\tilde{g}+\bigg(-\frac{\partial}{\partial r}\left(\frac{\beta \rho_{0}}{\alpha+\beta g_{0}}\right)-\frac{\partial g_{0}}{\partial r}\cdot\frac{\rho_{0}\beta^2}{(\alpha+\beta g_{0})^2}
\end{equation}
$$-\frac{1}{r}\cdot\frac{\beta\rho_{0}}{\alpha+\beta g_{0}}\bigg)\frac{\partial \tilde{g}}{\partial r}-\left(\frac{\beta \rho_{0}}{\alpha+\beta g_{0}}\right)\frac{\partial^2 \tilde{g}}{\partial r^2}-\left(\frac{1}{r^2}\cdot\frac{\beta \rho_{0}}{\alpha+\beta g_{0}}\right)\frac{\partial^2 \tilde{g}}{\partial \theta^2}$$ \\
\noindent We now assume solutions of the form
$$\tilde{\rho}=e^{st}e^{in\theta}F(r),$$
$$\tilde{g}=e^{st}e^{in\theta}G(r),$$
$$\tilde{v}_{r}=e^{st}e^{in\theta}H_{r}(r),$$
\begin{equation}
\tilde{v}_{\theta}=e^{st}e^{in\theta}H_{\theta}(r).
\end{equation}
\noindent which we will use to prove asymptotic stability. The variables $s$ and $n$ are arbitrary parameters and $F$, $G$, $H_{r}$, and $H_{\theta}$ are functions of $r$. These solutions give rise to the following useful properties:

\begin{center}
$$\frac{\partial \tilde{\rho}}{\partial t} = s\tilde{\rho},\hspace{0.8mm} \frac{\partial \tilde{\rho}}{\partial \theta} = in\tilde{\rho}, \hspace{0.8mm}  \frac{\partial^2 \tilde{\rho}}{\partial \theta^2} = -n^2 \rho $$
\vspace{-0.1cm}
$$\frac{\partial \tilde{g}}{\partial t} = s\tilde{g}, \hspace{0.8mm}  \frac{\partial \tilde{g}}{\partial \theta} = in\tilde{g}, \hspace{0.8mm} \frac{\partial^2 \tilde{g}}{\partial \theta^2} = -n^2 g$$
\vspace{-0.1cm}
$$\frac{\partial \tilde{v}_{r}}{\partial t} = s\tilde{v}_{r}, \frac{\partial \tilde{v}_{r}}{\partial \theta} = in\tilde{v}_{r}, \frac{\partial^2 \tilde{v}_{r}}{\partial \theta^2} = -n^2 \tilde{v}_{r},$$
\begin{equation}
\frac{\partial \tilde{v}_{\theta}}{\partial t} = s\tilde{v}_{\theta}, \frac{\partial \tilde{v}_{\theta}}{\partial \theta} = in\tilde{v}_{\theta}, \frac{\partial^2 \tilde{v}_{\theta}}{\partial \theta^2} = -n^2 \tilde{v}_{\theta}.
\end{equation} \\
\end{center} 

\noindent Plugging in these properties into our equation, we get \\
$$s\tilde{\rho}+\tilde{v}_{r}\frac{\partial \rho_{0}}{\partial r}+(inv_{0,\theta})\tilde{\rho}=\bigg(\frac{\partial}{\partial r}\left(\frac{\beta}{\alpha+\beta g_{0}}\right)\cdot\frac{\partial g_{0}}{\partial r}-\left(\frac{\beta}{\alpha+\beta g_{0}}\right)\bigg(\frac{1}{r}\frac{\partial g_{0}}{\partial r}$$
$$+\frac{\partial^2 g_{0}}{\partial r^2}\bigg)\bigg)\tilde{\rho}+\bigg(\frac{1}{r}+\frac{\partial g_{0}}{\partial r}\cdot\frac{\partial}{\partial r}\left(\frac{\beta}{\alpha+\beta g_{0}}\right)\bigg)\frac{\partial \tilde{\rho}}{\partial r}+\frac{\partial^2 \tilde{\rho}}{\partial r^2}+\bigg(-\frac{\partial g_{0}}{\partial r}$$
\begin{equation}
\cdot\frac{\partial}{\partial r}\left(\frac{\rho_{0}\beta^2}{(\alpha+\beta g_{0})^2}\right)+\frac{\beta^2\rho_{0}}{\alpha+\beta g_{0}}\cdot\left(\frac{1}{r}\frac{\partial g_{0}}{\partial r}+\frac{\partial^2 g_{0}}{\partial r^2}\right)-\frac{n^2}{r^2}\cdot\frac{\beta \rho_{0}}{\alpha+\beta g_{0}}\bigg)\tilde{g}
\end{equation}
$$-\left(\frac{\partial}{\partial r}\left(\frac{\beta \rho_{0}}{\alpha+\beta g_{0}}\right)+\frac{\partial g_{0}}{\partial r}\cdot\frac{\rho_{0}^2 \beta}{(\alpha+\beta g_{0})^2}+\frac{1}{r}\frac{\beta \rho_{0}}{\alpha+\beta g_{0}}\right)\frac{\partial \tilde{g}}{\partial r}-\left(\frac{\beta \rho_{0}}{\alpha+\beta g_{0}}\right)\frac{\partial^2 \tilde{g}}{\partial r^2},$$ \\
\noindent which we can rearrange with a little algebra to obtain \\
$$\frac{\partial \tilde{\rho}}{\partial t}+\tilde{v}_{r}\frac{\partial \rho_{0}}{\partial r}+v_{0,\theta}\frac{\partial \tilde{\rho}}{\partial \theta}=\frac{1}{r}\frac{\partial \tilde{\rho}}{\partial r}+\frac{\partial^2 \tilde{\rho}}{\partial r^2}-\frac{\partial}{\partial \theta}\left(\frac{\beta \rho_{0}}{\alpha+\beta g_{0}}\right)\cdot\frac{\partial \tilde{g}}{\partial r}+\frac{\partial g_{0}}{\partial r}\bigg(\frac{\beta}{\alpha+\beta g_{0}}\frac{\partial \tilde{\rho}}{\partial r}$$
\begin{equation}
+\frac{\partial}{\partial r}\left(\frac{\beta}{\alpha+\beta g_{0}}\right)\cdot\tilde{\rho}-\frac{\rho_{0}\beta^2}{(\alpha+\beta g_{0})^2}\frac{\partial \tilde{g}}{\partial r}-\frac{\partial}{\partial r}\left(\frac{\rho_{0}\beta^2}{(\alpha+\beta g_{0})^2}\right)\tilde{g}\bigg)-\left(\frac{\beta \rho_{0}}{\alpha+\beta g_{0}}\right)\bigg(\frac{1}{r}\frac{\partial \tilde{g}}{\partial r}
\end{equation}
$$+\frac{\partial^2 \tilde{g}}{\partial r^2}+\frac{1}{r^2}\frac{\partial^2 \tilde{g}}{\partial \theta^2}\bigg)-\left(\frac{\beta \tilde{\rho}}{\alpha+\beta g_{0}}-\frac{\beta^2\rho_{0}\tilde{g}}{(\alpha+\beta g_{0})^2}\right)\left(\frac{1}{r}\frac{\partial g_{0}}{\partial r}+\frac{\partial^2 g_{0}}{\partial r^2}\right).$$ \\
\noindent Finally, we plug the relation $\tilde{g}=\left(\frac{\lambda}{s+1}\right)\tilde{\rho}$ from Section \ref{continuum} (equation \ref{eq:61}) and rearrange to get an equation only in terms of $\tilde{\rho}$ and its $r$ derivatives:

$$s\tilde{\rho}=\left(\frac{\partial}{\partial r}\left(\frac{\beta}{\alpha+\beta g_{0}}\right)\cdot\frac{\partial g_{0}}{\partial r}-\left(\frac{\beta}{\alpha+\beta g_{0}}\right)\left(\frac{1}{r}\frac{\partial g_{0}}{\partial r}+\frac{\partial^2 g_{0}}{\partial r^2}\right)-inv_{0,\theta}\right)\tilde{\rho}+\bigg(\frac{1}{r}$$
$$+\frac{\partial g_{0}}{\partial r}\cdot\frac{\partial}{\partial r}\left(\frac{\beta}{\alpha+\beta g_{0}}\right)\bigg)\frac{\partial \tilde{\rho}}{\partial r}+\frac{\partial^2 \tilde{\rho}}{\partial r^2}+\bigg(-\frac{\partial g_{0}}{\partial r}\cdot\frac{\partial}{\partial r}\left(\frac{\rho_{0}\beta^2}{(\alpha+\beta g_{0})^2}\right)+\frac{\beta^2\rho_{0}}{\alpha+\beta g_{0}}\cdot\bigg(\frac{1}{r}\frac{\partial g_{0}}{\partial r}$$

$$+\frac{\partial^2 g_{0}}{\partial r^2}\bigg)-\frac{n^2}{r^2}\cdot\frac{\beta \rho_{0}}{\alpha+\beta g_{0}}\bigg)\tilde{g}-\left(\frac{\partial}{\partial r}\left(\frac{\beta \rho_{0}}{\alpha+\beta g_{0}}\right)+\frac{\partial g_{0}}{\partial r}\cdot\frac{\rho_{0}^2 \beta}{(\alpha+\beta g_{0})^2}+\frac{1}{r}\frac{\beta \rho_{0}}{\alpha+\beta g_{0}}\right)\frac{\partial \tilde{g}}{\partial r}$$
\begin{equation}
-\left(\frac{\beta \rho_{0}}{\alpha+\beta g_{0}}\right)\frac{\partial^2 \tilde{g}}{\partial r^2}-\left(\frac{\partial \rho_{0}}{\partial r}\right)\tilde{v}_{r}.
\end{equation}

\subsection{Velocity Equation, $\vec{v}$}
We start with the velocity equation $\vec{v}$ defined in Section \ref{continuum}: \\
\begin{equation}
\frac{d\vec{v}}{dt}=\frac{\partial \vec{v}}{\partial t}+(\vec{v} \cdot \nabla)\vec{v}=b \nabla g. \\
\end{equation}
\noindent Taking the material derivative $(\vec{v}\cdot\nabla)\vec{v}$, we obtain \\
\begin{equation}
\frac{\partial \vec{v}}{\partial t}+\left(v_{r}\frac{\partial v_r}{\partial r}+\frac{v_{\theta}}{r} \frac{\partial v_r}{\partial \theta}-\frac{v_{r} v_{\theta}}{r}\right)\hat{r}+\left(v_{r}\frac{\partial v_{\theta}}{\partial r}+\frac{v_{\theta}}{r}\frac{\partial v_{\theta}}{\partial r}+\frac{v_{r}v_{\theta}}{r}\right)\hat{\theta}=b \nabla g, \\
\end{equation}
\noindent which we rearrange so that all $\hat{r}$ and $\hat{\theta}$ terms are grouped to get \\
$$\left(\frac{\partial v_{r}}{\partial t}+v_{r}\frac{\partial v_r}{\partial r}+\frac{v_{\theta}}{r} \frac{\partial v_r}{\partial \theta}-\frac{v_{r} v_{\theta}}{r}\right)\hat{r}+\left(\frac{\partial v_{\theta}}{\partial t}+v_{r}\frac{\partial v_{\theta}}{\partial r}+\frac{v_{\theta}}{r}\frac{\partial v_{\theta}}{\partial r}+\frac{v_{r}v_{\theta}}{r}\right)\hat{\theta}$$
\begin{equation}
=\left(b\frac{\partial g}{\partial r}\right)\hat{r}+\left(\frac{b}{r}\frac{\partial g}{\partial \theta}\right)\hat{\theta}. \\
\end{equation}
\noindent Separating the $\hat{r}$ and $\hat{\theta}$ components, we have the following two equations: \\
\begin{equation} \label{eq:88}
\frac{\partial v_{r}}{\partial t}+v_{r}\frac{\partial v_r}{\partial r}+\frac{v_{\theta}}{r} \frac{\partial v_r}{\partial \theta}-\frac{v_{r} v_{\theta}}{r}=b\frac{\partial g}{\partial r},
\end{equation}
\begin{equation} \label{eq:89}
\frac{\partial v_{\theta}}{\partial t}+v_{r}\frac{\partial v_{\theta}}{\partial r}+\frac{v_{\theta}}{r}\frac{\partial v_{\theta}}{\partial r}+\frac{v_{r}v_{\theta}}{r}=\frac{b}{r}\frac{\partial g}{\partial \theta}. 
\end{equation} \\
\noindent We begin by linearizing equation \ref{eq:88}.

\noindent Substituting $v_{r} = v_{0,r}+\epsilon_{1}\tilde{v}_{r}$ and $v_{\theta} = v_{0,\theta}+\epsilon_{1}\tilde{v}_{\theta}$, we obtain \\
$$\frac{\partial v_{0,r}}{\partial t}+\epsilon_{1}\frac{\partial \tilde{v}_{r}}{\partial t}+(v_{0,r}+\epsilon_{1}\tilde{v}_{r})\left(\frac{\partial v_{0,r}}{\partial r}+\epsilon_{1}\frac{\tilde{v}_{r}}{\partial r}\right)$$
\begin{equation}
\cdot\left(\frac{\partial v_{0,r}}{\partial \theta}+\epsilon_{1}\frac{\partial \tilde{v}_{r}}{\partial \theta}\right)-\frac{(v_{0,r}+\epsilon_{1}\tilde{v_{r}})(v_{0,\theta}+\epsilon_{1}\tilde{v}_{\theta})}{r}=b\left(\frac{\partial g_{0}}{\partial r}+\epsilon_{1}\frac{\partial \tilde{g}}{\partial r}\right)
\end{equation} \\
\noindent Noting that $v_{0,r}=0$ and eliminating all terms that are zeroth or second order in $\epsilon_{1}$ (since $\frac{\partial \vec{v}_{0}}{\partial t}=0$ and $\epsilon_{1}$ is small), we obtain \\
\begin{equation}
\epsilon_{1}\frac{\partial \tilde{v}_{r}}{\partial t}+\frac{\epsilon_{1}v_{0,\theta}}{r}\frac{\partial \tilde{v}_{r}}{\partial \theta}-\frac{2\epsilon_{1}v_{0,\theta}}{r}=\epsilon_{1}b\frac{\partial \tilde{g}}{\partial r}. 
\end{equation} \\
\noindent Dividing out $\epsilon_{1}$, we get \\
\begin{equation}
\frac{\partial \tilde{v}_{r}}{\partial t}+\frac{v_{0,\theta}}{r}\frac{\partial \tilde{v}_{r}}{\partial \theta}-\frac{2v_{0,\theta}}{r}=b\frac{\partial \tilde{g}}{\partial r}.
\end{equation}
\noindent We now substitute the solutions for $\tilde{v}_{r}$ and $
\tilde{v}_{\theta}$ in $(1)$ into the equation to obtain \\
\begin{equation}
s\tilde{v}_{r}+\frac{v_{0,\theta}}{r}\cdot in\tilde{v}_{r}-\frac{2v_{0,\theta}}{r}\tilde{v}_{\theta}=b\frac{\partial \tilde{g}}{\partial r},
\end{equation} \\
\noindent which we rearrange to get \\
\begin{equation}
\left(s+\frac{v_{0,\theta}}{r}\cdot in\right)\tilde{v}_{r}-\left(\frac{2v_{0,\theta}}{r}\right)\tilde{v}_{\theta}=b\frac{\partial \tilde{g}}{\partial r}.
\end{equation} \\
\noindent We now perform the same manipulations on equation \ref{eq:89}.
\noindent Substituting $v_{r} = v_{0,r}+\epsilon_{1}\tilde{v}_{r}$ and $v_{\theta} = v_{0,\theta}+\epsilon_{1}\tilde{v}_{\theta}$, we obtain
$$ \frac{\partial v_{0,\theta}}{\partial t}+\epsilon_{1}\frac{\partial \tilde{v}_{\theta}}{\partial t}+(v_{0,r}+\epsilon_{1}\tilde{v}_{r})\left(\frac{\partial v_{0,\theta}}{\partial r}+\epsilon_{1}\frac{\partial \tilde{v}_{\theta}}{\partial r}\right)+\left(\frac{v_{0,\theta}+\epsilon_{1}\tilde{v}_{\theta}}{r}\right)\left(\frac{\partial v_{0,\theta}}{\partial \theta}+\epsilon_{1}\frac{\partial \tilde{v}_{\theta}}{\partial \theta}\right)$$
\begin{equation}
+\frac{(v_{0,r}+\epsilon_{1}\tilde{v}_{r})(v_{0,\theta}+\epsilon_{1}\tilde{v}_{\theta})}{r}=\frac{b}{r}\left(\frac{\partial g_{0}}{\partial \theta}+\epsilon_{1}\frac{\partial \tilde{g}}{\partial r}\right). 
\end{equation} \\
\noindent Since $v_{0,r}=0$ and we eliminate all terms that are zeroth or second order in $\epsilon_{1}$, we get \\
\begin{equation}
\epsilon_{1}\frac{\partial \tilde{v}_{\theta}}{\partial t}+\epsilon_{1}\tilde{v}_{r}\frac{\partial v_{0,\theta}}{\partial r}+\frac{\epsilon_{1}v_{0,\theta}}{r}\frac{\partial \tilde{v}_{\theta}}{\partial \theta}+\frac{\epsilon_{1} v_{0,\theta}}{r}\tilde{v}_{r}=\epsilon_{1}\frac{b}{r}\frac{\partial \tilde{g}}{\partial \theta}.
\end{equation} \\
\noindent Dividing out $\epsilon_{1}$ gives us \\
\begin{equation}
\frac{\partial \tilde{v}_{\theta}}{\partial t}+\tilde{v}_{r}\frac{\partial v_{0,\theta}}{\partial r}+\frac{v_{0,\theta}}{r}\frac{\partial \tilde{v}_{\theta}}{\partial \theta}+\frac{v_{0,\theta}}{r}\tilde{v}_{r}=b\frac{\partial \tilde{g}}{\partial \theta}
\end{equation} \\
\noindent Finally, plugging in the solutions for $\tilde{v}_{r}$ and $\tilde{v}_{\theta}$ as above and rearranging gives us \\
$$s\tilde{v}_{\theta}+\tilde{v}_{r}\frac{\partial v_{0,\theta}}{\partial r}+\frac{v_{0,\theta}}{r}(in\tilde{v}_{\theta})+\frac{v_{0,\theta}}{r}\tilde{v}_{r}=in\frac{b}{r}\tilde{g}$$
\begin{equation}
\left(\frac{\partial v_{0,\theta}}{\partial r}+\frac{v_{0,\theta}}{r}\right)\tilde{v}_{r}+\left(s+\frac{inv_{0,\theta}}{r}\right)\tilde{v}_{\theta} = in\frac{b}{r}\tilde{g}.
\end{equation} \\
We therefore have the following system of linear equations in $\tilde{v}_{r}$, $\tilde{v}_{\theta}$, and $\tilde{g}$:
\begin{equation}
\left(s+\frac{v_{0,\theta}}{r}\cdot in\right)\tilde{v}_{r}-\left(\frac{2v_{0,\theta}}{r}\right)\tilde{v}_{\theta}=b\frac{\partial \tilde{g}}{\partial r}
\end{equation}
\begin{equation}
\left(\frac{\partial v_{0,\theta}}{\partial r}+\frac{v_{0,\theta}}{r}\right)\tilde{v}_{r}+\left(s+\frac{inv_{0,\theta}}{r}\right)\tilde{v}_{\theta} = in\frac{b}{r}\tilde{g}.
\end{equation} \\
\noindent We solve this system and plug in the relation $\tilde{g}=\left(\frac{\lambda}{s+1}\right)\tilde{\rho}$ (equation \ref{eq:61}) to get \\ 
\begin{equation} 
\tilde{v}_{\theta}=\left(\frac{\lambda}{s+1}\right)\frac{\left(\left(\frac{bin\tilde{\rho}}{r}\right)\left(s+\frac{inv_{0,\theta}}{r}\right)-\left(b\frac{\partial \tilde{\rho}}{\partial r}\right)\left(\frac{\partial v_{0,\theta}}{\partial r}+\frac{v_{0,\theta}}{r}\right)\right)}{\left(\frac{2v_{0,\theta}}{r}\right)\left(\frac{\partial v_{0,\theta}}{\partial r}+\frac{v_{0,\theta}}{r}+\left(s+\frac{inv_{0,\theta}}{r}\right)^2 \right)},
\end{equation}
\begin{equation}
\tilde{v}_{r}=\left(\frac{\lambda}{s+1}\right)\frac{\left(bs\frac{\partial \tilde{\rho}}{\partial r}+\left(\frac{binv_{0,\theta}}{r}\right)\left(\frac{\partial \tilde{\rho}}{\partial r}+\frac{2\tilde{\rho}}{r}\right)\right)}{\left(\frac{2v_{0,\theta}}{r}\right)\left(\frac{\partial v_{0,\theta}}{\partial r}+\frac{v_{0,\theta}}{r}+\left(s+\frac{inv_{0,\theta}}{r}\right)^2 \right)}.
\end{equation} \\
which we used in Section $3.3$.

% might want to show substitution of expression for tilde{v}_{\theta} back into Equation 2, and then say you can solve for tilde{v}_{r} in terms of only \tilde{g}... i.e. break up the third line of text from the bottom.

\end{document}